\documentclass{article} %

\usepackage{url} 

\usepackage{prooftree} 
\usepackage{gb4e} 

\usepackage{url} 

\title{Some  proof theoretical remarks\\ on quantification in ordinary language} 
\author{Michele Abrusci\thanks{Dipartimento di Filosofia, Universit{\`a} di Roma Tre} \and Christian Retor{\'e}\thanks{Universit\'e de Bordeaux \& IRIT-C.N.R.S. \& LaBRI-C.N.R.S.}}

\renewcommand\l\lambda 
\newcommand\existe\exists
\newcommand\qqst\forall
\newcommand\implique\supset
\newcommand\et\land
\newcommand\ou\lor

\newcommand\seq{\mathrel{\vdash}}

\newcommand{\editout}[1]{}

\usepackage{latexsym} 
\usepackage{amssymb} 
\newcommand{\ltyn}{\ensuremath{\Lambda\mathsf{Ty}_n}}
\newcommand\TYN\ltyn

\newcommand\us{$^{\mathrm{(us)}}$}

\begin{document}
\maketitle 
\titlepage

\begin{center}\bf 
Some  proof theoretical remarks

on quantification in ordinary language
\end{center} 

\bigskip\bigskip 

\paragraph{Abstract:} This paper surveys the common approach  to quantification and generalised quantification in formal linguistics and philosophy of language. 
We point out how this general setting departs from empirical linguistic data, and give some hints for  a different view based on proof theory, which on many aspects gets closer to the language itself. 
We stress the importance of Hilbert's operator epsilon and tau for, respectively, existential and universal quantifications. Indeed, these operators help a lot to construct semantic representation close to natural language, in particular with quantified noun phrases as individual terms. We also define guidelines for the design of the proof rules corresponding to generalised quantifiers. 

\paragraph{R{\'e}sum{\'e}} Cet article dresse un rapide panorama de l'approche commune de la quantification g{\'e}n{\'e}ralis{\'e}e ou non en linguistique formelle et en philosophie du langage. Nous montrons que ce cadre g{\'e}n{\'e}ral est va parfois {\`a} l'encontre des donn{\'e}es linguistiques, et nous donnons quelques indications  pour une approche diff{\'e}rente bas{\'e}e sur la th{\'e}orie de la d{\'e}monstration, qui sur bien des points s'avère plus proche de la langue. Nous soulignons l'importance des op{\'e}rateurs tau et  epsilon de Hilbert qui rendent respectivement compte de la qualification universelle et existentielle.  En effet ces op{\'e}rateurs permettent de construire des des repr{\'e}sentations s{\'e}mantiques qui suivent la lange avec, en particulier des groupes nominaux quantifi{\'e}es qui soient des termes individuels. Nous donnons aussi des principes pour d{\'e}finir des règles de d{\'e}duction qui correspondent aux quantificateurs g{\'e}n{\'e}ralis{\'e}s.

\newpage 

\section{Foreword: empirical data on quantification}

Although quantification in logic, in linguistics and in philosophy of language are much studied  we think that the dominant model theoretic approach partly leaves out linguistic reality. 
\cite{PetersWesterstahl2006quantifiers,Steedman2012scope,szabolcsi2010quantification,gabbay2009quantification}

Quantifiers are quite common in ordinary language and even more in speciality languages like sciences, in particular mathematics... unfortunately for the linguists, the more difficult quantifiers are the most common in natural language. 
Indeed, as we shall see besides the classical quantifiers of mathematics \emph{for all $\forall$} and \emph{there exists 
$\exists$}, natural language makes use of many other quantifiers, some of which are only implicitly stated. Consider for instance the plural definite determiner \emph{the} it can be used as universal quantification as in \emph{The children fell asleep}, and it can be used as \emph{most} 
was in \emph{The Brits love France}\footnote{Unless otherwise stated our examples are from the web. When they are ours, e.g. for phrasing differently an example, we use the mark \us.}  which does not mean that \emph{the Brits} all do, but more likely that a large part of them do, or even that a short majority does as in \emph{France may have elected Fran\c{c}ois Hollande, but she really needs a Napoleon Bonaparte}. 
Let us first say a linguistically oriented word on the usual quantifiers, namely \emph{for all} and \emph{there exists}. 

The existential quantifier is omnipresent in natural language. As soon as one speaks about something, this "something" is existentially quantified as in \emph{A man enters. He whistles.}\us The introduction of new referents, that are existentially quantified, and that are the beginning of anaphoric chains are an important par of the structure of a discourse. For instance, observe that the scope of the existential quantifier extends to the next sentence --- this has been a reason to for introducing various dynamic logics. The behaviour of existential quantifiers actually lead to a formalism which is structured by existential quantification, namely Discourse Representation Theory, which construct logical formula in a way that better matches the linear progression of sentences. It is often opposed to compositional semantics, since DRT proceeds top down from larger units to the smallest, words, but there by now exist compositional formulations of DRT. So existential quantifiers are 
really present in discourse and conversation, and they contribute to the structure and coherence of a discourse. 
\cite{KR93}

The universal quantifier, as it is use in mathematics,  is rather rare in ordinary conversation, as corpora study shows, an exception being  the reference to a group that can be inferred from the context  like in \emph{They all went to bed.}. The negative formulation is more frequent \emph{no one, nothing} etc. and if \emph{no one came} it means that \emph{every one did not come}.  

When speaking about natural language it is often assumed that the individuals that is being referred to by quantification is finite. \footnote{Another problem might be that the quantity a quantifier refers to can be  continuous, as in \emph{We drank all the wine.}}  This is not the case. Firstly we meets such expression in mathematical discussions, but also various cases: 

\begin{exe} 
\ex \emph{He wrapped up by explaining the dark future for the Universe when all the stars go away.} 
\ex \emph{All atoms are made from the same bits, which are called subatomic particles} 
\ex \emph{Just about all sentences in the English language fall into ten patterns determined by the presence and functions of nouns, verbs, adjectives, and adverbs.} 
\ex \emph{All ideas are welcome.}
\end{exe}

The most that we discussed above is also quite common and also applies to infinite collection of objects: 

\begin{exe} \it 
\ex In basic math, we're taught multiplication tables. We learn that most numbers are the answer to at least two different multiplication problems, some numbers are the answer to several, and then...
\ex Any module of known $\beta$ is weak. Most numbers have even $\beta$ and most of them are not antisymmetric.
\ex The number one reason why most people fail at dieting is simple: they can't stick to it.
\ex Since most numbers are not prime, it would waste time to check every number.
\ex ... thus, in the limit most numbers are not prime.
\end{exe} 

The human processing of quantifiers is quite difficult, in particular their alternation. There are experiments 
comforting this viewpoint, see e.g.  \cite{Szymanik2010cog}. 
But, on the other hand, one finds rarely but sometimes Henkin quantifiers in natural language, like 
\emph{Every member of the lab knows a member of every village sports club.}
What do people mean when using such quantifiers that would require a lifetime verification? This is a mystery, we know  that our study, as well as others,  via a confronting the formal model with what we hear or with psycholinguistics 
experiments may give some hints on the complexity of the human processing. 

These few introductory word should have convinced the reader that  quantification is an important and common phenomenon in natural langue, although mathematical logic ignores it (an important exception being Hilbert, as we shall see), and although natural language processing consider it as rare and irrelevant (an exception being categorial grammar). So we agree with semanticists: quantification is an important question connecting logic, linguistics and philosophy. 

\section{Remarks on the two mathematised quantifiers, \emph{there exists} and \emph{for all}} 

The two quantifiers that are used in mathematics, namely $\forall$ and $\exists$ on individuals 
received a proper logical description,
both in syntactical terms with rules (since Aristotle) and in semantical terms with interpretations in models (since the beginning of the XXth century). An important result is the coincidence of the two notions for first order logic known as the completeness theorem, established by G\"odel in 1930. It says that a formula is  true in every model (1) if and only if it is provable (2) --- the part (2) implies (1), known as soundness is easier to proved by induction on the derivation. 
Other famous theorems are the compactness and Löwenheim-Skolem theorems. 
The first one says that if any finite subset of a set $F$ of formulae is contradiction free (admits a model), so is $F$,
and the second ones say that a theory with an infinite model and a countable language there are models of every cardinality larger than countable. \cite{KK86}

Those results are well known --- especially completeness --- and there is by now not much proof theoretical study of quantification. 
But one should be aware that this correspondence only works for first-order classical logic --- and not for higher order logic. For instance when second or higher order quantification is used, sub-Boolean algebras have to be considered. When other logics are considered trickier structures are required to maintain this correspondence. Even in the well studied case of intuitionistic logic, much more complicated models have to used, e.g. (pre)sheaves of L-structures or Kripke models. 
\cite{gabbay2009quantification}

Even if we only consider first order quantification, the correspondence between provability and truth is not as transparent as one may hope. Natural language provides good insights to perceive the distinction between the proof theoretical and the model theoretical approach. More formally stated, completeness expresses the coincidence between 
\begin{itemize} 
\item provability $\seq \forall x.\ P(x)$ inferred from a proof without a free hypothesis  involving a free $x$ (a proof with a \emph{generic} $x$ (the precise rule is given thereafter)  
\item 
and truth in every model, that is $\&_{x\in I}P(x)$ in every domain $I$
\end{itemize} 
An important remark is that the proof theoretical viewpoint is finite, while the second one involves two infinite sets: models are usually  infinite\footnote{In first order logic, one can say that the model has less than $n$ elements, but there is not way to say that the model is finite since it requires second order logic}, and the family of all models is usually an infinite set as well. 
The proof theoretic approach makes use of a generic element. To prove something about any number, 
the proof starts with \emph{Let $n$ be any number\ldots} which means that $n$ as no other property than being a number ---this methods goes back at least to Pythagora, (-580 – -490), long before any logical formalisation. 

It should be observed that natural language makes a distinction between the two notion of quantification above: \emph{each} and \emph{all} rather concern the complete enumeration of the elements in a collection while \emph{every}, \emph{any} or bare plurals (e.g. \emph{Ducks lay eggs}) rather concern the generic elements, law and universal rules. 
In natural language quantification is commonly formulated with a \emph{restriction} to a given class: \emph{All the stars go away},
and this is also the case in usual mathematics: $$\forall x\in \mathbb{N}\ \exists u_1 \in \mathbb{N}\ \exists u_2 \in \mathbb{N}\ \exists u_3 \in \mathbb{N}\ \exists u_4 \in \mathbb{N} \ x=(u_1)^2(u_1)^2(u_2)^2(u_3)^2(u_4)^2$$
When the restriction is a restriction to a well identified class like (inanimate) things or human beings, 
there are quantifiers without restriction to specific classes like \emph{everyone} or \emph{everything} are possible too: they  apply to a single property. 

The restriction to a given class in quantification disappeared from mathematical logic (but not from type theory) because of the following observation of Frege. Indeed, for philosophical reasons, he insisted in having a single one-sorted universe, and  using implicitly a correspondence between a set $M$ the corresponding predicate ($M(x)$), one is able to represent quantification with restriction to a given class: 

$$\forall x\in M\ P(x) \quad \equiv \forall x\ (M(x)\Rightarrow  P(x))$$

$$\exists x\in M\ P(x) \quad \equiv \exists x\ (M(x) \& P(x))$$

To present examples from the most studied quantifiers, we follow Arsitotle square of oppositions: 
\begin{itemize} 
\item All As are Bs
\item Some As are Bs 
\item No As are Bs
\item Not alls A are Bs
\end{itemize} 
In the original formulation, the expression $A$ and $B$ denote \emph{terms}, which are much vaguer than properties or predicates --- in the Middle-Age, the theory of \emph{suppositiones} already contributed to make the notion of \emph{term} neater. 

Existential quantification (Some As are Bs) is formulated with \emph{some}, \emph{there is}, \emph{a} it can also be formulated without restriction to a class of individuals by \emph{someone}, or \emph{something}.  The last two statements exhibit that natural language is finer grained, and that it makes a difference because of the focus between statement that usual logic consider as equivalent.

\begin{exe} 
\ex \emph{There's a tramp sittin' on my doorstep}
\ex \emph{Some girls give me money}
\ex \emph{Something happened to me yesterday} 
\ex \emph{Some politicians are crooks}
\ex \emph{? Some crooks are politicians}\us 
\end{exe}

Observe that the existential negative  (the fourth corner of Aristotle square of oppositions) is NOT lexicalised (in English and other languages). Psycholinguists know that they  is harder to grasp, rather ambiguous,  and in ordinary conversation an existential negative statement  is often understood wrongly as universal negative statements (e.g. example \ref{notevpic}  is rather often understood as \ref{nopic}). 
Also observe in our unambiguous rephrasing  \ref{somopnot} of \ref{notevop}, the focus has changed from the ones 
that can be funded to the ones that cannot, making the rephrased statement tougher.

\begin{exe}
\ex \label{notevpic} 
\emph{Not Every Picture Tells a Story} 
\ex \label{nopic} 
\emph{No picture tells a story}\us 
\ex 
\emph{Some Students Do Not Participate In Group Experiments Or Projects.} 
\ex \label{notevop}
\emph{Everyone is entitled to an opinion, but not every opinion is entitled to student government funding.}
\ex \label{somopnot} 
\emph{Everyone is entitled to an opinion, but some opinions are not entitled to student government funding.}\us 
\end{exe}

We already provided some examples of universal statements, here are some more, ranging on a potentially infinite domain:

\begin{exe}
\ex 
\emph{Each star in the sky is an enormous glowing ball of gas.} 
\item
\emph{All groups of stars are held together by gravitational forces.}
\ex 
\emph{Terence Tao, a Fields medalist, has published a paper that proves that every odd number greater than 1 is the sum of at most five primes.}
\end{exe} 

The universal negative statements are expressed either by \emph{no}, and without restriction to a class  by \emph{no one}, \emph{nothing}, … Here are some examples:  

\begin{exe} 
\ex \emph{Because no planet's orbit is perfectly circular, the distance of each varies over the course of its year.}
\ex \emph{Nothing's gonna change my world.} 
\ex \emph{Porterfield went where no colleague had gone previously this season, realising three figures.}  
\end{exe}

From a proof theoretical viewpoint how do $\exists$ and $\forall$ work? They both have introduction and elimination rules: 

\begin{itemize} 
\item Universal quantifier
\begin{itemize}
\item The $\forall$ introduction rule says as above that when a property has been established for an $x$ which does not enjoy any particular property (i.e. is not free in any hypothesis), one can conclude that the property holds for all individuals: 
$$\begin{prooftree} 
\begin{array}[b]{c}
\mbox{no free occurrence of $x$}\\ 
\mbox{in any $H_i$}\\ 
H_1,\ldots,H_n\seq P(x)   
\end{array} 
\justifies 
H_1,\ldots,H_n\seq \forall x.\ P(x) 
\using \forall_i
\end{prooftree}$$ 
\item The $\forall$ elimination rule says that when a property has been established for all individuals it can be inferred from any particular individual:  
$$\begin{prooftree} 
H_1,\ldots,H_n\seq \forall x.\ P(x) 
\justifies 
H_1,\ldots,H_n\seq P(a) 
\using \forall_e
\end{prooftree}$$ 
\end{itemize} 
\item Existential quantifier
\begin{itemize}
\item The $\exists$ introduction rule says that whenever a property $P$ 
holds of an individual $a$ one can infer that there exists an individual enjoying $P$ 
$$\begin{prooftree} 
H_1,\ldots,H_n\seq P(a)
\justifies 
H_1,\ldots,H_n\seq \exists x.\ P(x) 
\using\exists_i
\end{prooftree}$$ 
\item The $\exists$ elimination rule is trickier.  It says that if assuming $P(x)$ and nothing more about $x$ we derived $C$ which does not depend on $x$, and if we have a proof of $\exists x.\ P(x)$ we can conclude, without the hypothesis $P(x)$ that $C$: 
$$\begin{prooftree} 
H'_1,\ldots,H'_p\seq \exists x.\ P(x) 
\qquad 
\begin{array}[b]{c}
\mbox{no free occurrence of $x$}\\
\mbox{in any $H_i$ nor in $C$}\\ 
P(x),H_1,\ldots,H_n\seq C 
\end{array}  
\justifies 
H'_1,\ldots,H'_p, H_1,\ldots,H_n\seq C
\using \exists_e
\end{prooftree}$$ 
\end{itemize} 
\end{itemize} 

An alternative way to deal with quantification is to use individual concepts that are properties holding of exactly one individual. This is the reason why in Discourse Representation Theory some people write $John(x)$ while others write $\mathbf{j}$: the former notation describes an individual concept $John(\_)$ which has the property of being true of exactly one individual, and this can be stated in second and first order logic: $X$ is an individual concept 
whenever $C(X)=(\forall x \forall y (X(x) \land X(y) \Rightarrow  x = y)) \land \exists x.\ X(x)$ --- there are variants in which the existence is not demanded. Universal quantification  can be phrased with individual concepts via second order quantification. Indeed,   $\forall x.\ P(x)$ can be expressed as 
$\forall X (C(X)\Rightarrow P^\#(X))$ with $P^\#(X)= \exists X.\ (X(x) \& P(x))$, and, conversely, given a property $Q$ of individual concepts, the quantification  $\forall X.\ (C(X)\Rightarrow Q(X))$ corresponds to $\forall x.\ Q^b(x)$ with $Q^b(x)=
\exists X.\ (C(X)\land X(x) \land Q(X))$. The equivalence are established by formal proof in second order logic, with second order quantification, and by duality the same result holds for existential quantification. There is a variant according to which the non emptiness is not required, since it is possible to name individuals that do not exist — in this case one from of quantification is stronger than the other. 
One may wonder why using individual concepts and therefore use  second order logic. A reason is that it is closer to the notion of term of Ancient and Medieval logic: both \emph{Socrates} and \emph{human beings} are \emph{terms}. A more concrete reason is that, in common possible worlds semantics with rigid designators, it is impossible to interpret \emph{I do not believe that Tully is Cicero}.   \cite{Lacroix2011master}

Proof theoretic approaches naturally lead to rules of refutation. 
How does one refute a universal statement $\forall x P(x)$? Of course it amounts to try to prove the negation, that is an existential statement $\exists x \lnot P(x)$, but there are different ways to do so. Consider the following 
example:

\begin{exe}
\ex \emph{[The AKC notes] that any dog may bite.} 
\ex \emph{No, Rex would never bite.}\us 
\ex \emph{Basset hounds do not bite.}\us 
\end{exe} 

The difference is that the first answer picks up an element in the relevant model, while the second answer remains with generic elements. This is related to the Avicenian idea that a property of a term (individual or not) is always asserted for the term as part of a class: it is more related to type theory than to the Fregean view of a single universe. 

\section{Hilbert's operators: $\iota,\epsilon,\tau,\ldots$} 

After the quantifier \emph{the one and unique individual such that $P$ \ldots} introduced by Russell for definite description,
Hilbert (with Ackerman and Bernays) intensively used \emph{generic} elements for quantification, the study of which culminated in the second volume of \emph{Grundlagen der Mathematik} \cite{HBvol2}.  It should be stressed that these operators are introduced and describe with natural language examples, which is not so common in Hilbert's writings. 
We shall first present the $\epsilon$ operator which recently lead to important work in linguistics in particular with von Heusinger work. \cite{EgliHeusinger1995,Heusinger1997,Heusinger2004}

The first step due to Russell was to denote by $\iota_x.\ F$ the unique individual enjoying the property $F$ 
in a definite description like the first sentence and to remain undetermined when existence and uniqueness do not hold. 
\cite{Russell1905}

\begin{exe} 
\ex \emph{\emph{The present president of France} was born in Rouen.}\us \glt (existence and uniqueness hold) 
\ex \emph{\emph{The present king of France} was born in Pau.}\us  \glt (existence fails) 
\ex \emph{\emph{The present minister} was born in Barcelona.}\us  \glt (uniqueness fails) 
\end{exe} 

From this idea, Hilbert introduced  an individual existential term 
defined from a formula: given a formula $F(x)$ with a free variable $x$ 
one defines the term $\epsilon_x.\ F$  in which the occurrences of $x$ in $F$ are bound 
(this is the original notation, nowadays this term is often written as $\epsilon_x.\ F$). 
Whenever an element, say $a$,  enjoys $F$, then the  epsilon term $\epsilon_x.\ F$ enjoys $F$. 

Dually, Hilbert introduced a universal generic element $\tau_x.\ F$, which corresponds to the generic elements used in mathematical proofs: if $\tau_x.\ F$ enjoys the property $F$ then every individual does.

The evident deduction rules for $\tau$ and $\epsilon$ are as follows: 
\begin{itemize} 
\item 
From $A(x)$ with $x$ generic in the proof (no free occurrence of $x$ in any hypothesis),  infer $A(\tau_x.\ A(x))$
\item 
From $B(c)$ infer $B(\epsilon_x.\ B(x))$.  
\end{itemize} 

The additional rules can be found by duality: 
\begin{itemize} 
\item 
From $A(x)$ with $x$ generic in the proof (no free occurrence of $x$ in any hypothesis),  infer $A(\epsilon_x.\ \lnot A(x))$
\item 
From $B(c)$ infer $B(\tau_x.\ \lnot B(x))$.  
\end{itemize} 

Hence we have $F(\tau_x.\ F(x))\equiv \forall x. F(x)$ and $F(\epsilon_x.\ F(x))\equiv \exists x.\ F(x)$ 
and because of negation, one only of these operator is needed, usually the $\epsilon$ operator and the logic is known as the epsilon calculus. .

Hilbert turned these symbols into a mathematically quite nice object, since it allows to fully 
describe quantification and with simple rules. The first and second epsilon theorem 
basically say that this is an alternative 
formulation of first order logic. 
\begin{description}
\item[First epsilon theorem] 
When inferring a formula $C$ without $\epsilon$ symbol nor quantifier from  formulae $\Gamma$ not involving the $\epsilon$ symbol nor quantifiers  the derivation can be done within quantifier free  predicate calculus. 
\item[Second epsilon theorem]
When inferring a formula $C$ without $\epsilon$ symbol from  formulae $\Gamma$ not involving the $\epsilon$ symbol  the derivation can be done within predicate calculus. 
\end{description}

This way, it provided the first correct proof of Herbrand theorem (much before mistakes where found and solved by Goldfarb) 
and a proof of the consistence of Peano arithmetic (at the same time as Gentzen did). 
The extension of the cut  elimination patterns to these operators does not seem too complicated 
and people already worked on them, by 
\cite{mints2008}

In Hilbert's book these operators are explained with natural language examples,
but a very important linguistic property is not stated: 
the $\epsilon_x F$ as a type (both in the intuitive and in the formal sense) of a noun phrase, 
and is meant to be the argument of a predicate (for instance the subject of a verb),
thus being a \emph{suppositio} in the medieval sense. 
\cite{libera1996querelle,KK86}

Nowadays there has been a renewed interest in these epsilon formulation of quantification,
in particular by von Heusinger. He uses a variant of the epsilon for definite descriptions, leaving out the uniqueness of the iota operator of Russell, one reason being that the context often determines a unique object, the most salient one. 
As say it is a variant because it is not clear whether one still has the equivalence with ordinary existential quantification: he constructs an epsilon term whenever there is an expression like \emph{a man} or \emph{the man}  but it is not clear how 
does one asserts that $man(\epsilon_x.\ man(x))$. The distinction between $\epsilon$ and $\eta$ is that the former selects the most salient, while the later selects a new one.

\section{Generalised quantifiers} 

The abundant literature on generalised quantifiers takes place in a rather conventional setting, along the lines drawn by in particular by Frege.  It is assumed that generalised quantifiers like \emph{thirty per cent} or \emph{many} or \emph{most} are functions of two predicates. Indeed, the aforementioned Fregean trick to handle restricted quantification in a one-sorted logic ($\forall x\in M\ P(x) \quad \equiv \forall x\ (M(x)\Rightarrow  P(x))$) 
does not apply: as an easy exercise shows, the two sentences below are not equivalent: 

\begin{exe}
\ex Most students go out on Thursday nights. 
\ex For most individuals, if they are students then they go out on Thursday nights. 
\end{exe}

They are classified according to their behaviour with respect to the two predicates: covariant, contravariant, and properties wrt. existential quantification (weakness). Some rules are even provided, but \emph{in fine} they amount to count how many elements are in the class corresponding to the restriction, in the main predicate itself and in their intersection. The proof rules that are given, at some point, always refers to the intended model, and at this point one can see that they are not proof rules in the ordinary sense. 
This sometimes leads to inaccurate results. For instance, \emph{most} is defined as \emph{the majority of} although their actual usage is quite different: one can assert the first sentence but not the second  one in:  

\begin{exe}
\ex 
\emph{The majority of French electors voted Hollande in 2012}
\ex 
\emph{Most French electors voted Hollande in 2012}
\end{exe} 

Indeed, \emph{most} is a vague quantifier, and as \emph{many} starts beyond a flexible and context dependent percentage. 
The fact that the cardinality approach, the usual one,  is wrong is proved by the following examples:

\begin{exe} 
\ex \emph{Most numbers are not prime.}
\ex \emph{Most people have children.} (no one can tell whether there is a finite number of them)
\end{exe} 

Indeed, what is clear is that a measure is need but counting is not the only way to measure a set. 
Otherwise the first sentence would be false, while it appears in an advanced book on number theory, 
and while anyone knowing some basic maths agrees with such a sentence. 

Another problematic point is that the usual vision of (generalised)  quantifiers is to consider them as functions of several predicates. 
The mathematical quantifier are functions of a \emph{single} predicate: given a predicate $P(x)$  one can form $\forall x.\ P(x)$. 
General quantifiers depends on several predicates, usuallly the restriction to a given class and the main predicate. 
The quantifier $Q$ (some, all, most,...)  in a sentence like \emph{Q children sleep} is supposed, in the common view, 
to be a function that applies to the predicates \emph{children} (restriction class) and \emph{sleep}. 
Hence \emph{Q children} is meaningless in this traditional view. 
This goes against the syntactic and cognitive structure of language. 
If a speaker starts a sentence with \emph{Q children}  before he came to the verb, the hearer already has in mind an image of  what is \emph{Q children} with respect to \emph{children}. 
This drawback of the common approach is avoided when using Hilbert operators or choice functions. 
Our argument can be considered as a formal developement of Geach ideas  on \emph{terms} and \emph{predicates} 
\cite{geach1962reference}

\section{Towards rules for some generalised quantifiers} 

With this idea in mind, we started to study quantification from two converging viewpoints. 
One viewpoint is the interface between syntax and semantics, 
i.e. how do we map the syntactic structure of a sentence 
to a logical formula that represents the semantics of the analysed sentence. 
The other viewpoint is how one handles such formulae: how do we prove and refute them, how do we interpret them? 

Regarding the first viewpoint, we advocate for a typed version of Hilbert's generic elements. 
Consider the following sentences: 
\begin{exe} 
\ex Most dogs bite. 
\ex Most of the students that passed algebra passed logic. 
\end{exe} 
In those sentences the first restriction class, namely \emph{dogs} can be viewed as a type, while in the second example the restriction class \emph{students that passed algebra} is a formula with a unique free variable. The operator we use is a constant in a typed lambda calculus, and actually we have two distinct versions of the operator introducing a \emph{most} generic element. A first one  maps a type to an element in this type (hence the type of the operator  is false, but there is no problem to have a constant with this type) and the second ones takes as its argument a predicate over the type \emph{student} and return a \emph{student} and adds a presupposition that this \emph{student} has \emph{passed algebra}. We propose to do so for all generalise quantifiers, that is to have generic elements corresponding to them, computed from typed Hilbert's style operators. 
This is a way to compute their semantic representations, but it does not say how they are interpreted (we only provided some hints), not how one argues with such quantifiers. 
\cite{Retore2012rlv}

On this later point, regarding the \emph{the majority of} quantifier, we thought about what could be rules and refutations for such a quantifier. There are at least two ways to refute that the majority of (meaning at least 50\%) 
the $A$ have the property $P$:
\begin{itemize}
\item 
Only the minority of the $A$ has the property $P$
\item 
There is another property $Q$ which holds for the majority of the $A$, with no $A$ satisfying both $P$ and $Q$. 
\end{itemize} 

To add a generalised quantifier to a proof system one should actually introduce a pair of dual quantifiers,
a variant $\forall^*$ of $\forall$  and a variant $\exists^*$ of $\exists$ 
Then one has to decide  one of the following two possibilities, the first one being unlikely:
\begin{itemize}
\item 
$\forall^*x.\ A(x)$ implies $\forall x.\ A(x)$ and so $\exists x.\ A(x)$  implies $\exists^*x.\ A(x)$ 
\item 
$\exists^*x.\ A(x)$ implies $\exists x.\ A(x)$ and so $\forall x.\ A(x)$  implies $\forall^*x.\ A(x)$ 
\end{itemize} 
In both the cases, one of the two new quantifiers is obtained by extending an existing rule, 
and the other one oby restricting an existing rule. Although there cannot be complete set of rules,
it is likely to find family of rules  and refutation techniques that would model a large part of those quantifier's behaviour. 

\section{Conclusion} 
After some critics of the common set theoretic approach, we gave some hints for a proof theoretical approach to quantification. One idea is to use typed versions of Hilbert operator for computing the semantic representations. 
The other idea is to use principles on proofs and cu-elimination to design the rules of generalised quantifiers. 

\bibliographystyle{plain}
\bibliography{bigbiblio} 

\end{document}